\documentclass[smallextended]{svjour3}       
\smartqed  
\usepackage{graphicx}
\begin{document}
\title{An efficient mathematically correct scale free CORDIC}

\author{Yassine HACHA\"ICHI \and Younes LAHBIB}

\institute{LAMSIN-ENIT \\
              Universit\'e de Tunis El Manar, TUNISIA.\\
              \email{Yassine.Hachaichi@ipeiem.rnu.tn}  \\
ENICarthage\\
Universit\'e de Carthage, TUNISIA.\\
\and Electronics and Micro-electronics Laboratory\\
Universit\'e de Monastir, TUNISIA.\\
\email{Younes.Lahbib@enicarthage.rnu.tn}}


\maketitle

\begin{abstract}
In order to approximate transandental functions, several algorithms were proposed.
Historically, polynomial interpolation, infinite series, $\cdots$ and other
$+,\times, -$ and $/$ based algorithms were studied for this purpose.

The CORDIC (COordinate Rotation DIgital Computer)
introduced by Jack E. Volder in 1959, and generalized 
by 
J. S. Walther a few years later, is a hardware based algorithm
for the approximation of trigonometric, hyperbolic and
logarithmic functions.

As a consequence, CORDIC is used for applications in
diverse areas such as signal and image processing.
For these reasons, several modified versions were proposed.

In this article, we present an
overview of the CORDIC algorithm for the computation of the circular
 functions, essentially the scaling free version,
and we will give a substential improvement to the commonly used one.
\end{abstract}
\section{Introduction}
In 1959, Volder \cite{VOL}, introduced the CORDIC algorithm in order
to compute approximations of trigonometric functions.
This method is still used because of its adequacy to hardware design.
It is a recursive method using only shift-and-add operations.

A decade later, J. S. Walther in \cite{WAL}, generalized this method to other
transendantal functions used in engineering fields.

The development of the CORDIC algorithm and architectures \cite{LD} has taken place for
achieving the highest throughput rate and reduction of hardware-complexity as well as
the computational latency of implementation. Some of the typical approaches for reducing
complexity implementation are targeted on minimization of using the scaling-operation and
complexity of barrel-shifters and adders in the CORDIC engine.
However, one of the problems
associated with the classical
CORDIC formulation is that the scale factor depends of the angle, and is not constant. The
complexity of the computation of the scale factor
is in principle comparable to that of the basic CORDIC
process itself.
In a recent work, \cite{GOR}, a new algorithm, CORDIC II, is proposed that substitutes the
CORDIC micro-rotation by a new angle set.

Aiming to eliminate scale multiplication in conventional CORDIC, 
scale free CORDIC was used to eliminate the scale factor, see the
piooneering papres \cite{AMK,AK} and also \cite{MRC,MTB,MBG}.
The scale free CORDIC algorithm for cosine and sine functions is proved to be
faster and efficient in terms of area and accuracy compared to conventional
CORDIC.

We give in this paper a method in order to minimize the number of iterations ine the
CORDIC method. This is given by computing the closest elementary angle to the
residual one at each iteration.
Our second contribution is the correction of the Taylor series used for the composed functions.
We will prove that with our polynomial approximation, we will get faster computation for the same acuity.
The CORDIC algorithm operates either in, rotation mode or vectoring
mode, following linear, circular or hyperbolic coordinate trajectories.
In this paper, we focus on rotation mode CORDIC using circular trajectories.

\section{The CORDIC algorithm.}
The idea behind conventional CORDIC algorithm is the rotation
of a vector $[x_{in}\ y_{in}]^T$ in cartesian coordinate which can be
expressed in (1), where $[x_{out}\ y_{out}]^T$ is the output vector 
produced after rotation and $\theta$ is the angle of rotation.
\begin{equation}
\left(\begin{array}{c}
x_{out}\\ y_{out}\\
\end{array}
\right)=
\left(\begin{array}{cc}
\cos(\theta)&-\sin(\theta)\\ \sin(\theta)&\cos(\theta)\\
\end{array}\right)
\left(\begin{array}{c}
x_{in}\\ y_{in}\\
\end{array}\right)
\end{equation}
This can also be written as
\begin{equation}
\left(\begin{array}{c}
x_{out}\\ y_{out}\\
\end{array}
\right)=\cos(\theta)
\left(\begin{array}{cc}
1&-\tan(\theta)\\ \tan(\theta)&1\\
\end{array}\right)
\left(\begin{array}{c}
x_{in}\\ y_{in}\\
\end{array}\right)
\end{equation}
We split the rotation angle in a
sum of angles, and carries out the rotation by a series of the
so called micro-rotation by these angles.
The idea is to decompose any angle $\theta$ into a sum of some ''elemntary'' angles
\begin{equation}
\theta=\alpha_1+\cdots +\alpha_n
\label{dec}
\end{equation}
where $\alpha_k=\pm\arctan(2^{-l})$.

If we use the fact that, if $R_\theta$ denotes the matrix of the 2D rotation of angle $\theta$:
$$R_\theta \times R_{\theta'}= R_{\theta+\theta'}$$
We can translate the equation (3) into the matrix product :
$$R_\theta=R_{\alpha_1}\times \cdots \times R_{\alpha_n}$$

\subsection{The conventional CORDIC}
The  conventional CORDIC method performs a sequence of rotations by elementary angles.
Any rotation $\theta$ on the plan can be decomposed into a composition (matrix product)
of $n$ elementary rotations.

When taking $\theta_k=\arctan(2^{-k})$, the equation (2) becomes:
$$R_{\theta_k}=\frac{1}{\sqrt{1+2^{2k}}}
\left(\begin{array}{cc}
1&-2^{-k}\\ 2^{-k}&1\\
\end{array}\right)$$
Using also the identity
$$R_{-\theta}=\cos(\theta)
\left(\begin{array}{cc}
1&\tan(\theta)\\ -\tan(\theta)&1\\
\end{array}\right)$$
we obtain, for $\varepsilon_k=\pm 1$,
$$R_{\varepsilon_k\theta_k}=\frac{1}{\sqrt{1+2^{2k}}}
\left(\begin{array}{cc}
1&-\varepsilon_k 2^{-k}\\ \varepsilon_k 2^{-k}&1\\
\end{array}\right)$$
The idea is that the angles used are constant, so we have a constant
scale $K=\prod_{k=1}^n\frac{1}{\sqrt{1+2^{2k}}}$, which approximately equals, according to the litterature \cite{SSG}, $0.60725$.
For this aim, we construct a sequence of vectors $[x_k\ y_k\ z_k]^T$ by the recurrence schema:
\begin{equation}
\left(\begin{array}{c}
x_0\\y_0\\z_0\\
\end{array}\right)=
\left(\begin{array}{c}
x_{in}\\y_{in}\\\theta\\
\end{array}\right)
\end{equation}
and
\begin{equation}
\left\{
\begin{array}{ccccc}
x_{k+1}&=&x_k&-&\varepsilon_k2^{-k}y_k\\
y_{k+1}&=&\varepsilon_k2^{-k}x_k&+&y_k\\
z_{k+1}&=&z_k&-&\varepsilon_k \theta_k
\end{array}
\right.
\end{equation}
After the fixed number of iterations, we mutiply the resulting vector by the constant $K$, this means
$[x_{out}\ y_{out}]=K [x_n\ y_n]$.

The essence of the CORDIC algorithm is that he is multiplication free (only shift-and-add operations). The scale multiplication,
also called compensation, to get an output vector isometric to the input one, causes a problem.

The introduction of the scale free CORDIC is then legitimated.
\subsection{The correct scale free CORDIC for sine and cosine}
The scale free CORDIC for circular functions is based on the Taylor series 
$$\sin(x)=\sum_{n=0}^\infty \frac{(-1)^n\cdot x^{2n+1}}{(2n+1)!}=x-\frac{x^3}{6}+\frac{x^5}{120}+\cdots $$
$$\cos(x)=\sum_{n=0}^\infty \frac{(-1)^n\cdot x^{2n}}{(2n)!}=1-\frac{x^2}{2}+\frac{x^4}{24}+\cdots$$
$$\arctan(x)=\sum_{n=0}^\infty \frac{(-1)^n\cdot x^{2n+1}}{2n+1}=x-\frac{x^3}{3}+\frac{x^5}{5}+\cdots $$
The approximation Taylor polynomial of the \emph{composed functions} to order 5:
$$\sin(\arctan(x))\approx\left(x-\frac{x^3}{3}+\frac{x^5}{5}\right)-\frac{1}{6}\left(x-\frac{x^3}{3}+\frac{x^5}{5}\right)^3+$$
$$+\frac{1}{120}\left(x-\frac{x^3}{3}+\frac{x^5}{5}\right)^5$$
$$\cos(\arctan(x))\approx 1-\frac{1}{2}\left(x-\frac{x^3}{3}+\frac{x^5}{5}\right)^2+\frac{1}{24}\left(x-\frac{x^3}{3}+\frac{x^5}{5}\right)^4$$
When we truncate the polynomials to the order 5, we obtain the right equations:
\begin{equation}
\cos(\arctan(x))\approx 1-\frac{1}{2}x^2+\frac{3}{8}x^4
\end{equation}
\begin{equation}
\sin(\arctan(x))\approx x-\frac{1}{2}x^3+\frac{3}{8}x^5
\end{equation}
We can observe that $\sin(\arctan(x))=x\cos(\arctan(x))$, this is simply due to the fact that 
\begin{equation}
\frac{\sin(\arctan(x))}{\cos(\arctan(x))}=\tan(\arctan(x))=x
\end{equation}
We the obtain, for the elementary angles $\theta_k=\arctan(2^{-k})$, and remarking that $\frac{3}{8}=\frac{1}{4}+\frac{1}{8}$
$$\sin(\arctan(2^{-k}))\approx 2^{-k}-2^{-3k-1}+2^{-5k-2}+2^{-5k-3}$$
$$\cos(\arctan(2^{-k}))\approx 1-2^{-2k-1}+2^{-4k-2}+2^{-4k-3}$$
The rotation matrix $M_{\theta_k}$ becomes:
\begin{equation}
\left(\begin{array}{cc}
1-2^{-2k-1}+3\cdot 2^{-4k-3}&-2^{-k}+2^{-3k-1}-3\cdot 2^{-5k-3}\\
2^{-k}-2^{-3k-1}+3\cdot 2^{-5k-3}&1-2^{-2k-1}+3\cdot 2^{-4k-3}
\end{array}
\right)
\end{equation}
As we know, all the works we have seen uses the Taylor series for sine and cosine functions and replace
$\theta_k=\arctan(2^{-k})$ by $2^{-k}$, see \cite{AMK,AK} and also \cite{MRC,MTB,MBG}. 
The error is that when using a Taylor polynomial of a composite function $f\circ g$,
we have to use the same degree and truncate the resulting polynomial at the demanded degree, you can see
\cite{MON}.

In order to give an empirical proof, we will compare the orders 3, 4 and 5 of ou method 
 to the recent works \cite{MRC,AK}.
\section{Benchmark of scale free CORDIC for circular functions}
In order to minimize the number of iterations of the CORDIC algorithm, we choose the microrotations to be the closest $\arctan(2^{-k})$ to the residual angle.
This can be done by choosing the closest power of 2 to $\tan(\theta)$, where {$\theta$} represent the risidual angle.

Due of the continuity of the function $\arctan$, if $\tan(\theta)$ is close to $2^{-k}$, then so is $\theta$ to $\theta_k$.
This leads us to choose $k$ such that, the closest $\theta_k =\arctan(2^{-k})$ to $\theta$  the following way:
\begin{equation}
k= Round\left(\log_2\left(\frac{1}{|\theta|}\right)\right)
\label{close}
\end{equation}
 we replace $\arctan(\theta)$ by $\theta$ without any loss of acuity
 because $\theta$ is very close to $\arctan(\theta)$ for $\theta$ in $[0,\frac{\pi}{4}]$.

As an example:
$$\frac{\pi}{16}= \theta_1 - \theta_4 - \theta_7 - \theta_{10} + \theta_{12} (\pm 10^{-5})$$
For a hardware design, the translation of our method for the binary representation
$|\theta|=0.\varepsilon_1 \cdots \varepsilon_i \cdots$ where $\varepsilon_i\in\{0,1\}$
\begin{itemize}
\item if for $i\geq 1$ we have $\forall j<i ; \varepsilon_j=0$ and $\varepsilon_i=1$ and $\varepsilon_{i+1}=0$ then $k=i$
\item if for $i\geq 1$ we have $\forall j<i ; \varepsilon_j=0$ and $\varepsilon_i=\varepsilon_{i+1}=1$ then $k=i-1$
\end{itemize}
In this section, we will compare our approximation to the one given in \cite{AK} and \cite{MRC} for order 3 Taylor approximations.
The range of angles used is $[0,\frac{\pi}{4}]$.
This range is enough, using the trigonometric identities, to can calculate any sine or cosine of any angle.

In \cite{AK}, $\cos(\arctan(2^{-k}))$ is approximated by $1-2^{-2k-1}$, and $\sin(\arctan(2^{-k}))$ is approximated by $2^{-k}-2^{-3k-3}$.
The authors in \cite{MRC}, use the approximation $\cos(\arctan(2^{-k}))\approx 1-2^{-2k-1}$, and $\sin(\arctan(2^{-k}))\approx 2^{-k}-2^{-3k-2}$.
The proposed method use the approximation $\cos(\arctan(2^{-k}))\approx 1-2^{-2k-1}$, and $\sin(\arctan(2^{-k}))\approx 2^{-k}-2^{-3k-1}$,
which is the mathematically correct developpement.

The quadratic errors for cosine and sine function for the three method are summarized in the tables (\ref{tab:cos} and \ref{tab:sin}).
\begin{figure}[h]
\setlength{\unitlength}{1cm}
\begin{picture}(7,6)(0,0)
\put(4,5){\includegraphics[scale=0.5]{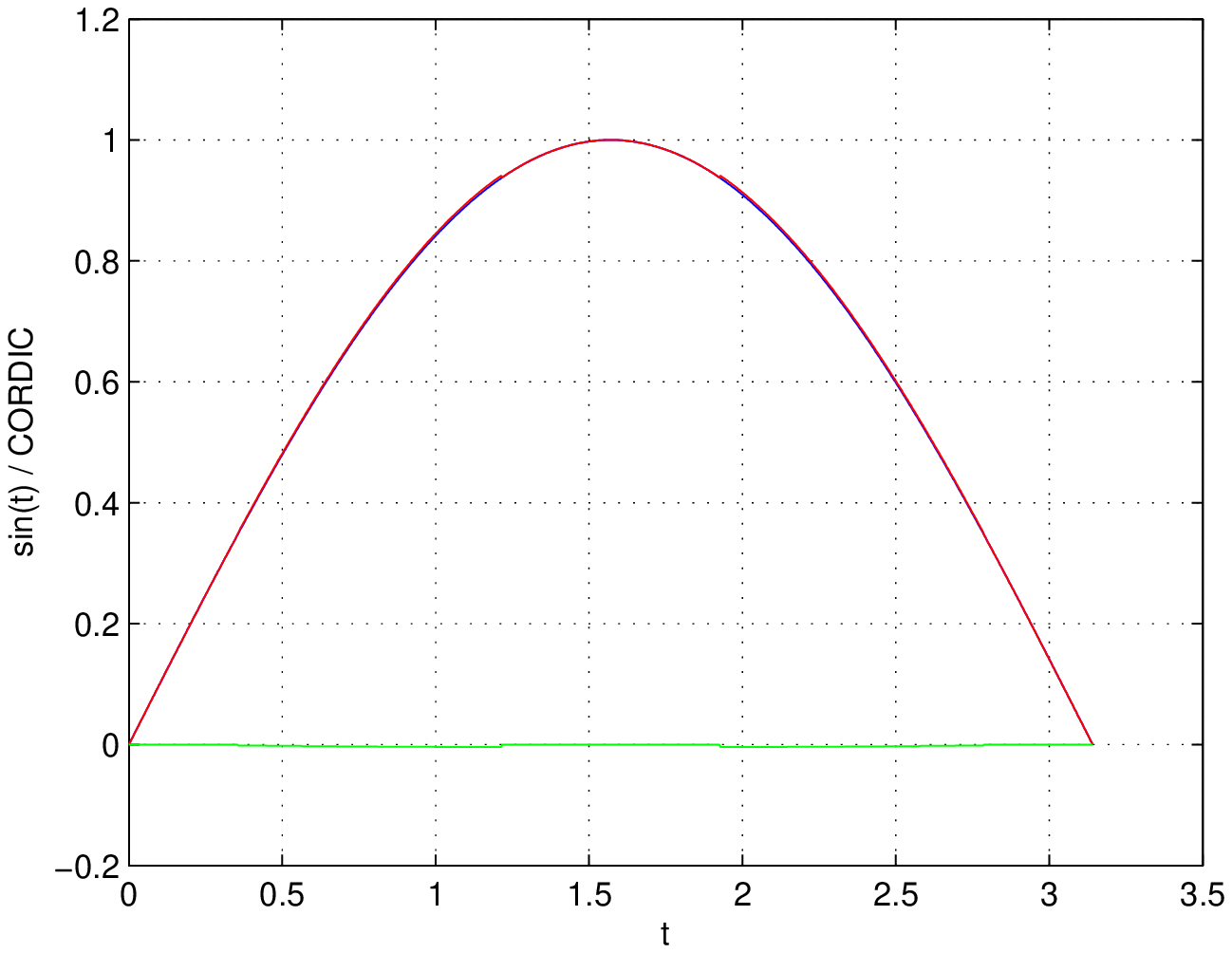}}
\end{picture}
\caption{\emph{$\sin$ approximation (in green the difference)}}
\label{fig:sin}
\end{figure}

\begin{table}[h]
\centering
\begin{tabular}{|c|c|c|c|}
\hline
&Method in \cite{AK} &Method in \cite{MRC}& Our method\\
\hline
3 iterations&3.0689e-004&2.9954e-004&2.8737e-004\\
\hline
4 iterations&2.3721e-004&2.2688e-004&2.0907e-004\\
\hline
5 iterations&2.2780e-004&2.1798e-004&2.0142e-004\\
\hline
\end{tabular}
\caption{Comparison of methods for cosine function}
\label{tab:cos}
\end{table}

\begin{figure}[h]
\setlength{\unitlength}{1cm}
\begin{picture}(7,6)(0,0)
\put(4,5){\includegraphics[scale=0.5]{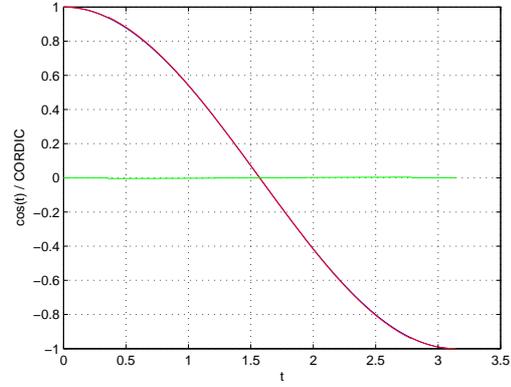}}
\end{picture}
\caption{\emph{$\cos$ approximation (in green the difference)}}
\label{fig:cos}
\end{figure}
\begin{table}[h]
\centering
\begin{tabular}{|c|c|c|c|}
\hline
&Method in \cite{AK} &Method in \cite{MRC}& Our method\\
\hline
3 iterations&0.001&8.0567e-004&4.9741e-004\\
\hline
4 iterations&8.2720e-004&5.4129e-004&5.3537e-005\\
\hline
5 iterations&8.0516e-004&5.1872e-004&5.5402e-005\\
\hline
\end{tabular}
\caption{Comparison of methods for sine function}
\label{tab:sin}
\end{table}

In the table (\ref{tab:DL}) below, we compare the quadratic errors of our method in different order Taylor approximation.
\begin{table}[!h]
\centering
\begin{tabular}{|cc|c|c|c|}
\hline
&&Order 5 &Order 4&Order 3\\
\hline
3 iterations &sin&5.0811e-004&4.9698e-004&4.9741e-004\\
&cos&1.6967e-004&1.6960e-004&2.8737e-004\\
\hline
4 iterations &sin&5.3826e-005&5.6214e-005&5.3537e-005\\
&cos&1.4850e-005&1.3429e-005&2.0907e-004\\
\hline
5 iterations &sin&8.4192e-006&5.0036e-005&5.5402e-005\\
&cos&1.0283e-005&1.1531e-005&2.0142e-004\\
\hline
\end{tabular}
\caption{Comparison of different orders of Taylor polynomials}
\label{tab:DL}
\end{table}

In figures (\ref{fig:cos} and \ref{fig:sin}), a MatLab simulation of our method is given.
In blue, the graph of our method, in red, the graph of the matlab function and in green the difference between them.
\section{Hardware implementation}
Common Hardware implementations of CORDIC algorithms are either iterative or pipelined \cite{YOU1,YOU2}. The main computation CORDIC unit is iterated in both cases. It is unrolled in the first class and rolled in the former using pipelined registers to store intermediate computations \cite{YOU3,YOU4}.   

A new design of the main computation unit is proposed in this paper and compared to the conventional CORDIC one. This is mainly a study to check if the theoretical results are feasible and simple to embed. Optimizations, complete CORDIC computation schemes, advanced CORDIC architectures and comparisons, which are based on the underlined computing unit, are future works.

The proposed scale-free CORDIC algorithm is based on Taylor polynomials. Three orders are evaluated for benchmarking the theoretical study (see Table \ref{tab:DL}). As the complexity of hardware architecture grows in function of the development order, only the order 3 is implemented in the hardware side. However, the impact of Taylor's order on hardware performances is also ongoing. It is not the focus of this paper.

\begin{figure}[!t]
\setlength{\unitlength}{1cm}
\begin{picture}(7,6)(0,0)
\put(2,0){\includegraphics[scale=0.5]{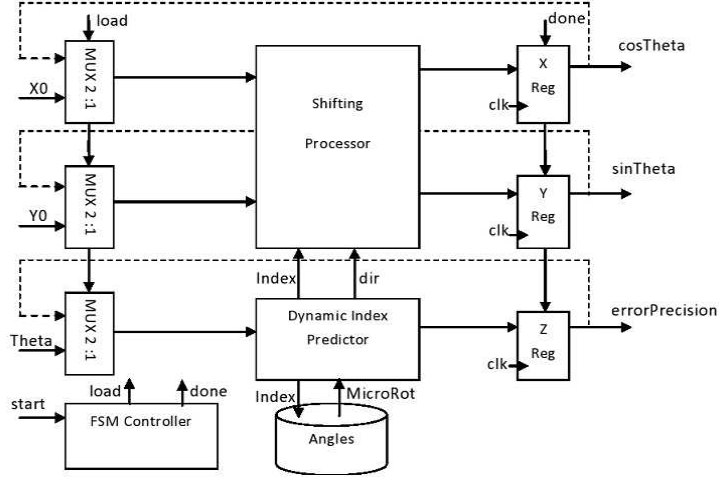}}
\end{picture}
\caption{\emph{Design of the CORDIC architecture.}}
\label{fig:arch}
\end{figure}

Hence, the implemented hardware architecture is restricted to the main computation unit and the order 3 of Taylor series. It is composed of 4 blocks: dynamic index predictor, shifting processor, storing angles ROM and FSM controller. Figure \ref{fig:arch} gives general description of the architecture. A detailed description is presented in the following sub-sections. The section ends with a summary of the main hardware results.
\subsection{FSM controller}
The controller is a finite state machine with three states. In the initial state, the signal {\em load} is set to initiate the initial values of the CORDIC core, namely $X_0=1, Y_0=0$ and the angle $\theta$. The second state is processed $2^N$ cycles according to an $N$-bits counter which fixes the number of CORDIC iterations.

For a given iteration, new intermediate values $X_n$ and $Y_n$ are obtained by shifting previous $X_{n-1}$ and $Y_{n-1}$ according to the closest
micro-rotation as given in equation \ref{close}.
 The third and final state sets the {\em done} signal. Cosine and sine of the angle $\theta$ are computed and stored in the output registers.     
\subsection{Dynamic Index predictor}
The main theoretical result proved in section 2. (see theoretical result) is implemented in the {\em Index Predictor}. The computation of the next index is the main improvement of the proposed hardware architecture. 
It estimates the optimal index with which we address the ROM and read the closest CORDIC micro-rotation for a given iteration. The determined angle is compared with the ongoing error $Z$ in order to compute the new direction of the micro-rotation. $X_n$ and $Y_n$ are then shifted index positions to right.
The next listing describes the behavior of this block. The hardware implementation is based on a 32 bits comparator. When the first sequence of {\tt '11'} bits are detected the block returns the corresponding $K$-index. Otherwise the first sequence of {\tt '10'} is cheeked, but the $(K-1)$-index is returned in this case. In this way, the most significant power of 2 micro-rotation is obtained.
{\begin{verbatim}
IndexEstim: process(Zerror)	
 begin
  for index in 31 downto 1 loop
    if(Zerror (index)='1' and Zerror (index-1) ='1') then
       Inxt <= std_logic_vector(to_unsigned(31-index,32));
       exit;
    elsif(Zerror(index) ='1' and Zerror(index-1) ='0') then
       Inxt<= std_logic_vector(to_unsigned(31-(index-1),4));
       exit;
    else
       Inxt <= (others => '0');
    end if;
   end loop;
end process IndexEstim;        
\end{verbatim}}
\subsection{Storing angles ROM}
Radix values of arctangent of $2^{-k}$ are stored in a ROM as constants, for $k$ within the range 0 to 31. The ROM values are not addressed inclemently as done in the conventional implementation of CORDIC. The closest 4 microrotations are read instead of 32 in the case of conventional CORDIC; The 4 indexes are estimated dynamically by the Index Predictor Block.

Listing below shows the addressing behavior the stored arctangent values.
\begin{verbatim}
IndexAccess: process(INDEX)
  begin
    case Index is	                           
       when "00000" => MicroRotI <= X"3243f6a8";
       _____

       when "11000" => MicroRotI <= X"0000003f";
       _____
       when "11111" => MicroRotI <= X"00000000";
    end Case;
  end process IndexAccess;
\end{verbatim}
\subsection{Shifting processor}
$X$ and $Y$ intermediate signals are computed as presented in listing below. This corresponds to processing equations in section 2. 
They are 32-bit coded signals in 2-complement format. Control signals {\em sgn} and {\em Index} are dynamically estimated and set by the 
Index predictor block. Note that listing below gives the behavior computing of the signal $Y$, it is similar in the case of $X$. 

The VHDL code behavior of the corresponding processing is similar to the listing below.
\begin{verbatim}
sinShift:process(Xsgn)
begin
  Yshift_tmp((31-I) downto 0) <= Xsgn(31 downto I);
  if Xsgn(31) = '0' then
	Yshift_tmp(31 downto (32-I)) <= (others => '0');
  else
    Yshift_tmp(31 downto (32-I)) <= (others => '1');
  end if;
end process sinShift;
\end{verbatim}

\subsection{Synthesis of the iterative CORDIC design}
The proposed iterative design based on the described computation unit is implemented on the Nexys3 spartan6 FPGA. 
Compared to the conventional architecture we reach the same cosine and sine results for only 4 iterations versus 32 iterations in the conventional case.
The precision error is about $10^{-3}$ in both cases. Hardware performances such as power, frequency and area are summarized in table \ref{tab:comp}.
\begin{table}[h]
\centering
\begin{tabular}{|c|c|}
\hline
Proposed Iterative CORDIC&Conventional Iterative CORDIC\\
\hline
\begin{tabular}{c|c|c}
Area & Power & Frequency \\
1020 Slices&100 mW&85.05 Mhz\\
\end{tabular}&
\begin{tabular}{c|c|c}
Area & Power & Frequency \\
541 Slices&72 mW&208.05 Mhz\\
\end{tabular}\\
\hline
\end{tabular}
\caption{Comparison of iterative architectures}
\label{tab:comp}
\end{table}
Implementation results show a degradation of the hardware performances of the proposed architecture. 
By an in-depth Analysis of sub-blocks, we find that the Index predictor block consumes alone 20 mw and uses 240 slices. 
The shifting processor block on the other side consumes also 20 mw and uses 540 slices.

The index predictor is an extra block in our case which explains the extra values against the conventional architecture. 
However, we think that the main reason is our coding style of the VHDL design which was behavioral. 
The behavioral synthesis infers the use of LUTs rather than basic logics. 
Hence, more optimized implementation should lead to less logic slices.
\subsection{Synthesis of the pipelined CORDIC design}
\begin{figure}[!t]
\setlength{\unitlength}{1cm}
\begin{picture}(7,6)(0,0)
\put(1,0){\includegraphics[scale=0.35]{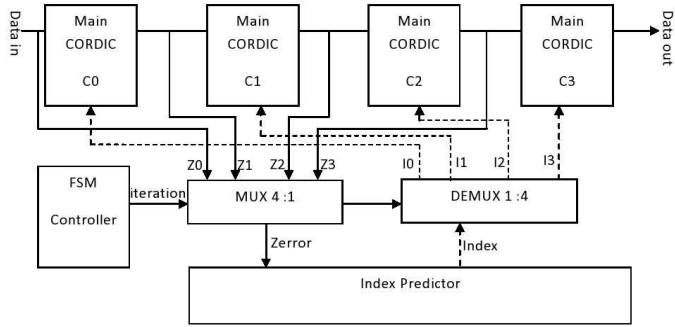}}
\end{picture}
\caption{\emph{Design of the pipelined CORDIC architecture.}}
\label{fig:pipe}
\end{figure}

A pipelined CORDIC architecture consists of rolling the main computation unit by storing intermediate computation into registers. 
In the case of the conventional architecture the main unit is rolled 32 times when processed data is 32-bit coded. 
The main improvement of our proposed architecture is rolling the same unit 3 or 4 times whatever the size of the processed data. 
More rolled units can be implemented if more precision is needed. With 3 units a precision of $10^{-2}$ is reached and $10^{-3}$ when 4 units are used. 
Figure \ref{fig:pipe} shows the proposed pipelined architecture. The index predictor which is resources consuming is instantiated only one time. 
The FSM controller enables the communication with only one pipeline stage. 
As shown in table \ref{tab:comp2} significant results compared to the iterative architecture are obtained. 
We save almost 50\% of area and power with a speedup of 10 Mhz, against only 16-bit conventional architecture.
\begin{table}[h]
\centering
\begin{tabular}{|c|c|}
\hline
Proposed Pipelined CORDIC&Conventional Pipelined CORDIC (16-BIT)\\
\hline
\begin{tabular}{c|c|c}
Area & Power & Frequency \\
4180 Slices&180 mW&55 Mhz\\
\end{tabular}&
\begin{tabular}{c|c|c}
Area & Power & Frequency \\
8200 Slices&300 mW&45 Mhz\\
\end{tabular}\\
\hline
\end{tabular}
\caption{Comparison of pipelined architectures}
\label{tab:comp2}
\end{table}
\section{Conclusion}
CORDIC algorithm has several applications in several domains, for an overview the reader can read \cite{50Y}.
The popularity of this method is due to the simplicity of its
hardware implementation, see \cite{LD} for example.

In this paper two improvements were made. First, we have minimized the number of iterations for some fixed error
by calculating the closest elementary angle to the residual one.
Second, we gave the correct polynomial approximation for the scale free CORDIC.
The comparison between our method and two other famous methods is given to
confirm empirically our theoretical proof.
In our simulation, we remark that the order of approximation of
Taylor series used meets the accuracy requirements.

In section 3, we showed that these methods have a simple hardware implementation, in order to meet the
objectives of the CORDIC's introduction.
The iterative and pipelined architectures were implemented, and significant improvements of hardware performance
were denoted in the pipelined case.
The future works will focus on the improvement of the hardware architecture.


\begin{thebibliography}{99}
\bibitem{AK}Aggarwal, S., \& Khare, K. (2012). \emph{Hardware efficient architecture for generating sine/cosine waves.} In VLSI Design (VLSID), 2012 25th International Conference on (pp. 57-61). IEEE.
\bibitem{AKHyp}Aggarwal, S., Meher, P. K., \& Khare, K. (2013). \emph{Scale-Free Hyperbolic CORDIC Processor and its
Application to Waveform Generation}. Circuits and Systems I: Regular Papers, IEEE Transactions on, 60(2), 314-326.
\bibitem{AMK}Aggarwal, S., Meher, P. K., \& Khare, K. (2012). \emph{Area-time efficient scaling-free cordic using generalized micro-rotation selection.} Very Large Scale Integration (VLSI) Systems, IEEE Transactions on, 20(8), 1542-1546.
\bibitem{YOU2} Antonius P. R., Nur A., Ashbir A. F., Naufal S., Trio A.,.(2015). \emph{FPGA implementation of CORDIC algorithms for sine and cosine generator.} IEEE Electrical Engineering and Informatics (ICEEI), Int. Conf. Proc. pp 1-6.
\bibitem{YOU3} Dian-Marie R., Scott M., Mihai S., Curran C., (2011), \emph{Design rules for implementing CORDIC on FPGAs.} IEEE Communications, Computers and Signal Processing (PacRim). Conf. Proc. pages 797 – 802.
\bibitem{GOR}Garrido, M., K\"allstr\"om, P., Kumm, M., \& Gustafsson, O., (2015). \emph{CORDIC II: A New Improved CORDIC Algorithm}. Circuits and Systems II: Regular Papers, IEEE Transactions on.
\bibitem{JSH}Jaime, F. J., S\`anchez, M. A., Hormigo, J., Villalba, J., \& Zapata, E. L. (2010). \emph{Enhanced scaling-free CORDIC.} Circuits and Systems I: Regular Papers, IEEE Transactions on, 57(7), 1654-1662.
\bibitem{LD}Lakshmi, B., \& Dhar,  A. S., (2010) \emph{CORDIC Architectures: A Survey,} Hindawi Publishing Corporation
VLSI Design, vol. 2010, pp. 1-19.
\bibitem{MBG} Maharatna, K., Banerjee, S., Grass, E., Krstic, M., \& Troya, A. (2005). Modified virtually scaling-free adaptive CORDIC rotator algorithm and architecture. Circuits and Systems for Video Technology, IEEE Transactions on, 15(11), 1463-1474.
\bibitem{MTB}Maharatna, K., Troya, A., Banerjee, S., \& Grass, E. (2004, November). \emph{Virtually scaling-free adaptive CORDIC rotator.}
 In Computers and Digital Techniques, IEE Proceedings- (Vol. 151, No. 6, pp. 448-456). IET.
\bibitem{50Y}Meher, P. K., Valls, J., Juang, T. B., Sridharan, K., \& Maharatna, K. (2009). \emph{50 years of CORDIC: Algorithms, architectures, and applications.} Circuits and Systems I: Regular Papers, IEEE Transactions on, 56(9), 1893-1907.
\bibitem{MRC}Mokhtar, A. S. N., Reaz, M. B. I., Chellappan, K., \& Ali, M. M. (2013). \emph{Scaling free CORDIC algorithm implementation of sine and cosine function.} In Proceedings of the World Congress on Engineering (WCE'13) (Vol. 2).
\bibitem{SSG}Surapong, P., Samman, F. A., \& Glesner, M. (2012). \emph{Design and Analysis of Extension-Rotation CORDIC Algorithms based on Non-Redundant Method.} International Journal of Signal Processing, Image Processing and Pattern Recognition, 5(1), 65-84.
\bibitem{MON}Monier, J. M., \emph{J'int\`egre. Analyse MPSI: cours et 1000 exercices corrig\'es.} Dunod Paris, 2003.
\bibitem{YOU4} Ray, A.,(1989). \emph{A survey of CORDIC algorithms for FPGA based computers.} International symposium on Field programmable gate arrays Pages 191-200.
\bibitem{YOU1}  Valls, J., Kuhlmann, M., Parhi, K. (2002). \emph{Evaluation of CORDIC Algorithms for FPGA Design.}
 Journal of VLSI signal processing systems for signal, image and video technology
November 2002, Volume 32, Issue 3, pp 207-222
\bibitem{VOL}Volder, J. E., (1959). \emph{The CORDIC trigonometric computing technique,} IRE
Trans. Electron. Comput., vol. EC-8, pp. 330–334.
\bibitem{WAL}Walther, J. S., (1971). \emph{A unified algorithm for elementary functions}, Spring
Joint Computer Conf., Proc., pp. 379-385
\end{thebibliography}
\end{document}